\documentclass[11pt, twoside]{article}
\usepackage{latexsym}
\usepackage{amsmath}
\usepackage{amssymb}
\usepackage[all]{xy}
\usepackage{amsfonts}
\usepackage{verbatim}
\usepackage{amsthm}
\usepackage{mathrsfs}
\usepackage{epsfig}
\usepackage{xy}
\usepackage{tensor}
\usepackage{array}
\usepackage{stmaryrd}
\usepackage{graphicx,color}
\usepackage{xcolor}
\usepackage[colorlinks=true,linkcolor=blue,citecolor=blue]{hyperref}
\usepackage{tikz}
\usetikzlibrary{arrows,calc}
\usepackage{etex}
\usepackage{mathdots}
\usepackage{float}
\usepackage{graphics}
\usepackage{pdflscape}
\usepackage{extarrows}
\usepackage{anysize,hyperref}
\input xypic
\xyoption{all}
\usepackage{bm}
\usepackage[perpage,symbol]{footmisc}
\usepackage{setspace}
\SelectTips{cm}{10}

\renewcommand{\cal}{\mathcal}
\def\A{\mathscr{A}}

\def\C{\mathscr{C}}

\def\Y{\mathscr{Y}}

\def\E{\mathbb{E}}
\def\s{\mathfrak{s}}
\def\Id{\mathrm{Id}}

\def\del{\delta}
\def\dr{\ar@{->}[r]}
\def\Im{\mbox{\rm Im}\,}\def\Ker{\mbox{\rm Ker}\,}

\def\Y{\mathscr{Y}}

\def\Hom{\mbox{Hom}}

\begin{document}
\baselineskip=15pt
\title{\Large{\bf  Torsion pairs and recollements of extriangulated categories\footnotetext{Yonggang Hu was supported by the National Natural Science Foundation of China (Grant  Nos. 11671126 and 12071120). Panyue Zhou was supported by the National Natural Science Foundation of China (Grant No. 11901190) and the Scientific Research Fund of Hunan Provincial Education Department (Grant No. 19B239).} }}
\medskip
\author{Jian He, Yonggang Hu and Panyue Zhou}

\date{}

\maketitle
\def\blue{\color{blue}}
\def\red{\color{red}}

\newtheorem{theorem}{Theorem}[section]
\newtheorem{lemma}[theorem]{Lemma}
\newtheorem{corollary}[theorem]{Corollary}
\newtheorem{proposition}[theorem]{Proposition}
\newtheorem{conjecture}{Conjecture}
\theoremstyle{definition}
\newtheorem{definition}[theorem]{Definition}
\newtheorem{question}[theorem]{Question}
\newtheorem{notation}[theorem]{Notation}
\newtheorem{remark}[theorem]{Remark}
\newtheorem{remark*}[]{Remark}
\newtheorem{example}[theorem]{Example}
\newtheorem{example*}[]{Example}

\newtheorem{construction}[theorem]{Construction}
\newtheorem{construction*}[]{Construction}

\newtheorem{assumption}[theorem]{Assumption}
\newtheorem{assumption*}[]{Assumption}

\baselineskip=17pt
\parindent=0.5cm

\begin{abstract}
\baselineskip=16pt
In this article, we prove that if  $(\mathcal A ,\mathcal B,\mathcal C)$ is  a recollement of extriangulated categories, then
torsion pairs in $\mathcal A$ and $\mathcal C$ can induce torsion pairs in $\mathcal B$, and the converse holds under natural assumptions. Besides,
we give mild conditions on a cluster tilting subcategory on the middle category of a recollement of extriangulated categories,  for the corresponding abelian quotients to form a recollement of abelian categories.\\[2mm]
\medskip
\end{abstract}
\pagestyle{myheadings}
\markboth{\rightline {\scriptsize J. He, Y. Hu and P. Zhou\hspace{2mm}}}
         {\leftline{\scriptsize  Torsion pairs and recollements of extriangulated categories}}

\section{Introduction}
The recollement of triangulated categories was introduced first by Beilinson, Bernstein, and Deligne, see \cite{BBD}. A fundamental example of a recollement situation of abelian categories appeared in the construction of perverse sheaves by MacPherson and Vilonen \cite{MV}. Recollements of  triangulated (abelian) categories can be viewed as `exact sequences' of triangulated (abelian) categories, which describe the middle term by a subcategory and a quotient category. It should be noted that these two recollement situations are now widely used in the study of  representation theory and algebraic geometry.

Let $(\mathcal A ,\mathcal B,\mathcal C)$ be a recollement of triangulated categories. Chen \cite{C} described how to glue together cotorsion pairs in $\mathcal A$ and $\mathcal C$ to obtain a cotorsion pair in $\mathcal B$.
Subsequently, Ma and Huang \cite{MH} did similar work  in the context of  abelian categories. More precisely,
 they showed how to construct a  torsion pair of the middle term  from that of two outer terms along  a recollement of abelian categories.

 Koenig and Zhu \cite{KZ} provided a general
framework for passing from triangulated categories to abelian categories by factoring out cluster tilting
subcategories.  It becomes a powerful tool to understand how  triangulated categories turn into abelian categories. Lin and Wang \cite{LW} used this method to construct  recollements of abelian categories from a recollement of triangulated categories.

The notion of extriangulated categories was introduced by Nakaoka and Palu in \cite{NP} as a simultaneous generalization of
exact categories and triangulated categories. Exact categories (abelian categories are also exact categories) and extension closed subcategories of an
extriangulated category are extriangulated categories, while there are some other examples of extriangulated categories which are neither exact nor triangulated, see \cite{NP,ZZ1,HZZ,NP1}. Hence many results hold on exact categories
and triangulated categories can be unified in the same framework.
Wang, Wei, and Zhang \cite{WWZ} introduced
the recollement of extriangulated categories, which is a simultaneous generalization of recollements of abelian categories and triangulated categories.
They also gave conditions such that the glued pair with respect to cotorsion pairs in
$\mathcal A$ and $\mathcal C$ is a cotorsion pair in $\mathcal B$ for a recollement $(\mathcal A ,\mathcal B,\mathcal C)$ of extriangulated categories. This result
 recovered a result given by Chen \cite{C} for the recollement of triangulated categories. But it cannot cover a result of Ma and Huang \cite{MH} for the recollement of abelian categories. Inspired by this, we consider to extend and generalize the related notions so as to cover some known results.

Let $\mathcal C$ be a triangulated category with a shift functor $[1]$. We need to pay attention to this fact: a pair $(\mathcal U,\mathcal V)$  of full subcategories of $\mathcal C$ is a cotorsion pair in \cite{N}
if and only if $(\mathcal U[-1], \mathcal V)$ is a torsion  pair in \cite{IY}.
But this fact is not necessarily true in an abelian category.
Based on this idea, we have a natural question of whether their results of Chen \cite{C} and Ma-Huang \cite{MH} can be unified under the framework of extriangulated categories. In this article, we give an affirmative answer.

Suppose that $\mathcal B$ admits a recollement relative to extriangulated categories $\mathcal A$ and $\mathcal C$.
Our first main result describes how to glue together torsion pairs $(\mathcal T_1,\mathcal F_1)$ in $\mathcal A$ and $(\mathcal T_2,\mathcal F_2)$ in $\mathcal C$, to obtain a torsion pair
$(\mathcal T,\mathcal F)$ of $\mathcal B$, see Theorem \ref{main}. This unifies their results of Chen \cite{C} and Ma-Huang \cite{MH} in the framework of extriangulated categories.
In the reverse direction, our second main result gives sufficient conditions on a torsion pair $(\mathcal T,\mathcal F)$ of $\mathcal B$, relative to the functors involved in the recollement, to induce torsion pairs in $\mathcal A$ and $\mathcal B$, see Theorem \ref{main1}.  Our third main result constructs a recollement
of abelian categories from a recollement of extriangulated categories, see Theorem \ref{zh}. This generalizes a result of Lin and Wang \cite{LW}.

This article is organized as follows. In Section 2, we give some terminologies
and some preliminary results. In Section 3, we prove our first and second main
results. In Section 4, we prove our third main
result. In Section 5, we give  an example to explain our main results.

\section{Preliminaries}
We briefly recall some definitions and basic properties of extriangulated categories from \cite{NP}.
We omit some details here, but the reader can find them in \cite{NP}.

Let $\mathcal{C}$ be an additive category equipped with an additive bifunctor
$$\mathbb{E}: \mathcal{C}^{\rm op}\times \mathcal{C}\rightarrow {\rm Ab},$$
where ${\rm Ab}$ is the category of abelian groups. For any objects $A, C\in\mathcal{C}$, an element $\delta\in \mathbb{E}(C,A)$ is called an $\mathbb{E}$-extension.
Let $\mathfrak{s}$ be a correspondence which associates an equivalence class $$\mathfrak{s}(\delta)=\xymatrix@C=0.8cm{[A\ar[r]^x
 &B\ar[r]^y&C]}$$ to any $\mathbb{E}$-extension $\delta\in\mathbb{E}(C, A)$. This $\mathfrak{s}$ is called a {\it realization} of $\mathbb{E}$, if it makes the diagrams in \cite[Definition 2.9]{NP} commutative.
 A triplet $(\mathcal{C}, \mathbb{E}, \mathfrak{s})$ is called an {\it extriangulated category} if it satisfies the following conditions.
\begin{itemize}
\item $\mathbb{E}\colon\mathcal{C}^{\rm op}\times \mathcal{C}\rightarrow \rm{Ab}$ is an additive bifunctor.

\item $\mathfrak{s}$ is an additive realization of $\mathbb{E}$.

\item $\mathbb{E}$ and $\mathfrak{s}$  satisfy the compatibility conditions in \cite[Definition 2.12]{NP}.
\end{itemize}

We collect the following terminology from \cite{NP}.

\begin{definition}
Let $(\mathcal{C},\E,\s)$ be an extriangulated category.
\begin{itemize}
\item[(1)] A sequence $A\xrightarrow{~x~}B\xrightarrow{~y~}C$ is called a {\it conflation} if it realizes some $\E$-extension $\del\in\E(C,A)$.
    In this case, $x$ is called an {\it inflation} and $y$ is called a {\it deflation}.

\item[(2)] If a conflation  $A\xrightarrow{~x~}B\xrightarrow{~y~}C$ realizes $\delta\in\mathbb{E}(C,A)$, we call the pair $( A\xrightarrow{~x~}B\xrightarrow{~y~}C,\delta)$ an {\it $\E$-triangle}, and write it in the following way.
$$A\overset{x}{\longrightarrow}B\overset{y}{\longrightarrow}C\overset{\delta}{\dashrightarrow}$$
We usually do not write this $``\delta"$ if it is not used in the argument.

\item[(3)] Let $A\overset{x}{\longrightarrow}B\overset{y}{\longrightarrow}C\overset{\delta}{\dashrightarrow}$ and $A^{\prime}\overset{x^{\prime}}{\longrightarrow}B^{\prime}\overset{y^{\prime}}{\longrightarrow}C^{\prime}\overset{\delta^{\prime}}{\dashrightarrow}$ be any pair of $\E$-triangles. If a triplet $(a,b,c)$ realizes $(a,c)\colon\delta\to\delta^{\prime}$, then we write it as
$$\xymatrix{
A \ar[r]^x \ar[d]^a & B\ar[r]^y \ar[d]^{b} & C\ar@{-->}[r]^{\del}\ar[d]^c&\\
A'\ar[r]^{x'} & B' \ar[r]^{y'} & C'\ar@{-->}[r]^{\del'} &}$$
and call $(a,b,c)$ a {\it morphism of $\E$-triangles}.

\item[(4)] An object $P\in\mathcal{C}$ is called {\it projective} if
for any $\E$-triangle $\xymatrix{A\ar[r]^{x}&B\ar[r]^{y}&C\ar@{-->}[r]^{\delta}&}$ and any morphism $c\in\mathcal{C}(P,C)$, there exists $b\in\mathcal{C}(P,B)$ satisfying $yb=c$.
We denote the full subcategory of projective objects in $\mathcal{C}$ by $\cal P$. Dually, the full subcategory of injective objects in $\mathcal{C}$ is denoted by $\cal I$.

\item[(5)] We say that $\mathcal{C}$ {\it has enough projectives}, if
for any object $C\in\mathcal{C}$, there exists an $\E$-triangle
$$\xymatrix{A\ar[r]^{x}&P\ar[r]^{y}&C\ar@{-->}[r]^{\delta}&}$$
satisfying $P\in\cal P$.  We can define the notion of having enough injectives dually.

\end{itemize}
\end{definition}

Assume that $(\mathcal C, \E, \s)$ is an extriangulated category.

\begin{proposition}\label{exact} \rm{\cite[Proposition 3.3]{NP} }
Let $\xymatrix{A\ar[r]&B\ar[r]&C\ar@{-->}[r]&}$
be an $\E$-triangle in $\mathcal C$. Then we have the following exact sequences:
$$\mathcal C(-, A)\xrightarrow{~}\mathcal C(-, B)\xrightarrow{~}\mathcal C(-, C)\xrightarrow{~}
\E(-, A)\xrightarrow{~}\E(-, B);
$$
$$\mathcal C(C,-)\xrightarrow{~}\mathcal C(B,-)\xrightarrow{~}\mathcal C(A,-)\xrightarrow{~}
\E(C,-)\xrightarrow{~}\E(B,-).$$
\end{proposition}

\begin{lemma}\label{00}{\rm \cite[Lemma 2.5]{WWZ}}
\rm Let $A\stackrel{f}{\longrightarrow}B\stackrel{g}{\longrightarrow}C\stackrel{}\dashrightarrow$ be an $\mathbb{E}$-triangle in $\mathcal{C}$. Then $f$ is an isomorphism if and only if $C\cong0$. Similarly, $g$ is an isomorphism if and only if $A\cong0$.

\end{lemma}

The following some notions of exact functors from \cite{WWZ}.

\begin{definition}\label{right} {\rm \cite[Definition 2.8]{WWZ}}
A morphism $f$ in $\mathcal{C}$ is called {\em compatible}, if ``$f$ is both an inflation and a deflation" implies that $f$ is an isomorphism.
\end{definition}

\begin{definition}\label{right}{\rm \cite[Definition 2.9]{WWZ}}
 A sequence $A\stackrel{f}{\longrightarrow}B\stackrel{g}{\longrightarrow}C$ in $\mathcal{C}$ is said to be {\em right exact} if
there exists an $\mathbb{E}$-triangle $K\stackrel{h_{2}}{\longrightarrow}B\stackrel{g}{\longrightarrow}C\stackrel{}\dashrightarrow$ and a deflation $h_{1}:A\rightarrow K$ which is compatible, such that $f=h_2h_1$. Dually one can also define the {\em left exact} sequences.

A $4$-term $\mathbb{E}$-triangle sequence $A{\stackrel{f}\longrightarrow}B\stackrel{g}{\longrightarrow}C\stackrel{h}{\longrightarrow}D$ is called {\em right exact} (resp. {\em left exact}) if there exist $\mathbb{E}$-triangles $A\stackrel{f}{\longrightarrow}B\stackrel{g_1}{\longrightarrow}K\stackrel{}\dashrightarrow$
and $K\stackrel{g_{2}}{\longrightarrow}C\stackrel{h}{\longrightarrow}D\stackrel{}\dashrightarrow$ such that $g=g_2g_1$ and $g_1$ (resp. $g_2$) is compatible.
\end{definition}
\begin{remark}\label{2.8} (1) A sequence $\eta: A\stackrel{f}{\longrightarrow}B\stackrel{g}{\longrightarrow}C$ is both left exact and right exact if and only  if $\eta$ is a conflation.

(2) If $\mathcal{C}$ is an abelian category, then $A\stackrel{f}{\longrightarrow}B\stackrel{g}{\longrightarrow}C$ is right exact if and only if $A\stackrel{f}{\longrightarrow}B\stackrel{g}{\longrightarrow}C{\longrightarrow}0$ is exact. Similarly, $A\stackrel{f}{\longrightarrow}B\stackrel{g}{\longrightarrow}C$ is left exact if and only if $0\stackrel{} {\longrightarrow} A\stackrel{f}{\longrightarrow}B\stackrel{g}{\longrightarrow}C$ is exact. If $\mathcal{C}$ is a triangulated category with the shift functor [1]. Then $\eta: A\stackrel{f}{\longrightarrow}B\stackrel{g}{\longrightarrow}C$  is right exact if and only if $A\stackrel{f}{\longrightarrow}B\stackrel{g}{\longrightarrow}C{\longrightarrow}A[1]$ is a triangle if and only if $\eta$ is left exact.
\end{remark}

\begin{definition}\label{right exact}{\rm \cite[Definition 2.12]{WWZ}} Let $(\mathcal{A},\mathbb{E}_{\mathcal{A}},\mathfrak{s}_{\mathcal{A}})$ and $(\mathcal{B},\mathbb{E}_{\mathcal{B}},\mathfrak{s}_{\mathcal{B}})$ be extriangulated categories. An additive covariant functor $F:\mathcal{A}\rightarrow \mathcal{B}$ is called a {\em right exact functor} if it satisfies the following conditions
\begin{itemize}
   \item [(1)] If $f$ is a compatible morphism in $\mathcal{A}$, then $Ff$ is compatible in $\mathcal{B}$.
   \item [(2)] If $A\stackrel{a}{\longrightarrow}B\stackrel{b}{\longrightarrow}C$ is right exact in $\mathcal{A}$, then $FA\stackrel{Fa}{\longrightarrow}FB\stackrel{Fb}{\longrightarrow}FC$ is right exact in $\mathcal{B}$. (Then for any $\mathbb{E}_{\mathcal{A}}$-triangle $A\stackrel{f}{\longrightarrow}B\stackrel{g}{\longrightarrow}C\stackrel{\delta}\dashrightarrow$, there exists an $\mathbb{E}_{\mathcal{B}}$-triangle $A'\stackrel{x}{\longrightarrow}FB\stackrel{Fg}{\longrightarrow}FC\stackrel{}\dashrightarrow$ such that $Ff=xy$ and $y: FA\rightarrow A'$ is a deflation and compatible. Moreover, $A'$ is uniquely determined up to isomorphism.)
 \item [(3)] There exists a natural transformation $$\eta=\{\eta_{(C,A)}:\mathbb{E}_{\mathcal{A}}(C,A)\longrightarrow\mathbb{E}_{\mathcal{B}}(F^{op}C,A')\}_{(C,A)\in{\mathcal{A}}^{\rm op}\times\mathcal{A}}$$ such that $\mathfrak{s}_{\mathcal{B}}(\eta_{(C,A)}(\delta))=[A'\stackrel{x}{\longrightarrow}FB\stackrel{Fg}{\longrightarrow}FC]$.

 \end{itemize}
Dually, we define the {\em left exact functor} between two extriangulated categories.
\end{definition}

\begin{definition}\label{exact functor}{\rm \cite[Definition 2.13]{WWZ}}
Let $(\mathcal{A},\mathbb{E}_{\mathcal{A}},\mathfrak{s}_{\mathcal{A}})$ and $(\mathcal{B},\mathbb{E}_{\mathcal{B}},\mathfrak{s}_{\mathcal{B}})$ be extriangulated categories.  We say an additive covariant functor $F:\mathcal{A}\rightarrow \mathcal{B}$ is an {\em exact functor} if the following conditions hold.
\begin{itemize}
  \item [(1)] If $f$ is a compatible morphism in $\mathcal{A}$, then $Ff$ is compatible in $\mathcal{B}$.
  \item [(2)] There exists a natural transformation $$\eta=\{\eta_{(C,A)}\}_{(C,A)\in{\mathcal{A}}^{\rm op}\times\mathcal{A}}:\mathbb{E}_{\mathcal{A}}(-,-)\Rightarrow\mathbb{E}_{\mathcal{B}}(F^{\rm op}-,F-).$$
  \item [(3)] If $\mathfrak{s}_{\mathcal{A}}(\delta)=[A\stackrel{x}{\longrightarrow}B\stackrel{y}{\longrightarrow}C]$, then $\mathfrak{s}_{\mathcal{B}}(\eta_{(C,A)}(\delta))=[F(A)\stackrel{F(x)}{\longrightarrow}F(B)\stackrel{F(y)}{\longrightarrow}F(C)]$.

  \end{itemize}
\end{definition}
\begin{proposition}\label{2.10} \rm{\cite[Proposition 2.14]{WWZ}} Let $(\mathcal{A},\mathbb{E}_{\mathcal{A}},\mathfrak{s}_{\mathcal{A}})$ and $(\mathcal{B},\mathbb{E}_{\mathcal{B}},\mathfrak{s}_{\mathcal{B}})$ be extriangulated categories. An additive covariant functor $F: \mathcal{A}\rightarrow\mathcal{B}$ is exact if and only if $F$ is both left exact and right exact.
\end{proposition}

\begin{remark}\label{111}If the categories $\mathcal{A}$ and $\mathcal{B}$ are abelian, Definition \ref{right exact} coincides with the usual right exact functor in abelian categories, and Definition \ref{exact functor} coincides with the usual exact functor. If the categories $\mathcal{A}$ and $\mathcal{B}$ are triangulated, by Remark \ref{2.8} and Proposition \ref{2.10}, we know that $F$ is a left exact functor if and only if $F$ is a triangle functor if and only if  $F$ is a right exact functor.
\end{remark}
Let us recall the concept of a recollement of extriangulated categories from \cite{WWZ}.

\begin{definition}\label{recollement}{\rm \cite[Definition 3.1]{WWZ}}
Let $\mathcal{A}$, $\mathcal{B}$ and $\mathcal{C}$ be three extriangulated categories. A \emph{recollement} of $\mathcal{B}$ relative to
$\mathcal{A}$ and $\mathcal{C}$, denoted by ($\mathcal{A}$, $\mathcal{B}$, $\mathcal{C}$), is a diagram
\begin{equation}\label{recolle}
  \xymatrix{\mathcal{A}\ar[rr]|{i_{*}}&&\ar@/_1pc/[ll]|{i^{*}}\ar@/^1pc/[ll]|{i^{!}}\mathcal{B}
\ar[rr]|{j^{\ast}}&&\ar@/_1pc/[ll]|{j_{!}}\ar@/^1pc/[ll]|{j_{\ast}}\mathcal{C}}
\end{equation}
given by two exact functors $i_{*},j^{\ast}$, two right exact functors $i^{\ast}$, $j_!$ and two left exact functors $i^{!}$, $j_\ast$, which satisfies the following conditions:
\begin{itemize}
  \item [(R1)] $(i^{*}, i_{\ast}, i^{!})$ and $(j_!, j^\ast, j_\ast)$ are adjoint triples.
  \item [(R2)] $\Im i_{\ast}=\Ker j^{\ast}$.
  \item [(R3)] $i_\ast$, $j_!$ and $j_\ast$ are fully faithful.
  \item [(R4)] For each $X\in\mathcal{B}$, there exists a left exact $\mathbb{E}_{\mathcal B}$-triangle sequence
  \begin{equation}\label{first}
  \xymatrix{i_\ast i^! X\ar[r]^-{\theta_X}&X\ar[r]^-{\vartheta_X}&j_\ast j^\ast X\ar[r]&i_\ast A}
   \end{equation}
  with $A\in \mathcal{A}$, where $\theta_X$ and  $\vartheta_X$ are given by the adjunction morphisms.
  \item [(R5)] For each $X\in\mathcal{B}$, there exists a right exact $\mathbb{E}_{\mathcal B}$-triangle sequence
  \begin{equation}\label{second}
  \xymatrix{i_\ast\ar[r] A' &j_! j^\ast X\ar[r]^-{\upsilon_X}&X\ar[r]^-{\nu_X}&i_\ast i^\ast X&}
   \end{equation}
 with $A'\in \mathcal{A}$, where $\upsilon_X$ and $\nu_X$ are given by the adjunction morphisms.
\end{itemize}
\end{definition}

\begin{remark}\label{77}(1) If the categories $\mathcal{A}$, $\mathcal{B}$ and $\mathcal{C}$ are abelian, then Definition \ref{recollement} coincides with the definition of recollement of abelian categories (cf. \cite{MH,P,FP}).

(2) If the categories $\mathcal{A}$, $\mathcal{B}$ and $\mathcal{C}$ are triangulated, then Definition \ref{recollement} coincides with the definition of recollement of triangulated categories (cf. \cite{BBD}).

(3) There exists an example of recollement of an extriangulated category which is neither abelian nor triangulated, for more details, see \cite{WWZ} and Example \ref{example1}.

\end{remark}
We collect some properties of a recollement of extriangulated categories, which will be used in the sequel.
\begin{lemma}\label{CY}\rm{\cite[Proposition 3.3]{WWZ}} Let ($\mathcal{A}$, $\mathcal{B}$, $\mathcal{C}$) be a recollement of extriangulated categories as \rm{(\ref{recolle})}.

$(1)$ All the natural transformations
$$i^{\ast}i_{\ast}\Rightarrow\Id_{\A},~\Id_{\A}\Rightarrow i^{!}i_{\ast},~\Id_{\C}\Rightarrow j^{\ast}j_{!},~j^{\ast}j_{\ast}\Rightarrow\Id_{\C}$$
are natural isomorphisms.

$(2)$ $i^{\ast}j_!=0$ and $i^{!}j_\ast=0$.

$(3)$ $i^{\ast}$ preserves projective objects and $i^{!}$ preserves injective objects.

$(3')$ $j_{!}$ preserves projective objects and $j_{\ast}$ preserves injective objects.

$(4)$ If $i^{!}$ (resp. $j_{\ast}$) is  exact, then $i_{\ast}$ (resp. $j^{\ast}$) preserves projective objects.

$(4')$ If $i^{\ast}$ (resp. $j_{!}$) is  exact, then $i_{\ast}$ (resp. $j^{\ast}$) preserves injective objects.

$(5)$ If $\mathcal{B}$ has enough projectives, then $\mathcal{A}$ has enough projectives; if $\mathcal{B}$ has enough injectives, then $\mathcal{A}$ has enough injectives .

$(6)$  If $\mathcal{B}$ has enough projectives and $j_{\ast}$ is exact, then $\mathcal{C}$ has enough projectives ; if $\mathcal{B}$ has enough injectives and $j_{!}$ is exact, then $\mathcal{C}$ has enough injectives.

$(7)$ If $\mathcal{B}$ has enough projectives and $i^{!}$ is  exact, then $\mathbb{E}_{\mathcal{B}}(i_{\ast}X,Y)\cong\mathbb{E}_{\mathcal{A}}(X,i^{!}Y)$ for any $X\in\mathcal{A}$ and $Y\in\mathcal{B}$.

$(7')$  If $\mathcal{C}$ has enough projectives and $j_{!}$ is  exact, then $\mathbb{E}_{\mathcal{B}}(j_{!}Z,Y)\cong\mathbb{E}_{\mathcal{C}}(Z,j^{\ast}Y)$ for any $Y\in\mathcal{B}$ and $Z\in\mathcal{C}$.

$(8)$  If $i^{\ast}$ is exact, then $j_{!}$ is  exact.

$(8')$ If $i^{!}$ is exact, then $j_{\ast}$ is exact.

\end{lemma}

\section{Glued torsion pairs}
\setcounter{equation}{0}

We first introduce the concept of torsion pairs in an extriangulated categories.
\begin{definition}\label{torsion}
Let $\mathcal{C}$ be an extriangulated category and $\mathcal{T}$, $\mathcal{F}\subseteq\mathcal{C}$ be a pair of subcategories of $\mathcal{C}$.
The pair $(\mathcal{T},\mathcal{F})$ is called a {\em torsion pair} in $\mathcal{C}$ if it satisfies the following conditions:

(a) $\Hom_\mathcal{C}(\mathcal{T},\mathcal{F})=0,$ that is, $\Hom_\mathcal{C}({T},{F})=0,$ for any $T\in\mathcal{T}$, $F\in\mathcal{F}$.

(b) For any $C\in\mathcal{C}$, there exists a conflation $T\longrightarrow C\longrightarrow F$ such that $T\in\mathcal{T}$ and $F\in\mathcal{F}$.
\end{definition}
\begin{remark}\label{222} (1) If the category $\mathcal{C}$ are abelian, then Definition \ref{torsion} coincides with the definition of torsion pair of abelian categories (cf. \cite{BBD,MH,D}).

(2) If the categoriy $\mathcal{C}$ are triangulated, then Definition \ref{torsion} coincides with the definition of torsion pair of  triangulated categories (cf. \cite{IY}).
\end{remark}

The following useful result can be found in \cite{WWZ}.

\begin{lemma}\label{1}\rm{\cite[Proposition 3.4]{WWZ}}
Let ($\mathcal{A}$, $\mathcal{B}$, $\mathcal{C}$) be a recollement of triangulated categories and $X\in\mathcal{B}$. Then the following statements hold.

$(1)$ If $i^{!}$ is exact, there exists an $\mathbb{E}_\mathcal{B}$-triangle
  \begin{equation*}\label{third}
  \xymatrix{i_\ast i^! X\ar[r]^-{\theta_X}&X\ar[r]^-{\vartheta_X}&j_\ast j^\ast X\ar@{-->}[r]&}
   \end{equation*}
 where $\theta_X$ and  $\vartheta_X$ are given by the adjunction morphisms.

$(2)$ If $i^{\ast}$ is exact, there exists an $\mathbb{E}_\mathcal{B}$-triangle
  \begin{equation*}\label{four}
  \xymatrix{ j_! j^\ast X\ar[r]^-{\upsilon_X}&X\ar[r]^-{\nu_X}&i_\ast i^\ast X \ar@{-->}[r]&}
   \end{equation*}
where $\upsilon_X$ and $\nu_X$ are given by the adjunction morphisms.
\end{lemma}

Our first main result is the following.

\begin{theorem}\label{main}
 Let $\mathcal A,\mathcal B,\mathcal C$ be three extriangulated category. Assume that
 $\mathcal B$ admits a recollement relative to $\mathcal A$ and $\mathcal C$ as follows
$$\xymatrix{\mathcal{A}\ar[rr]|{i_{*}}&&\ar@/_1pc/[ll]|{i^{*}}\ar@/^1pc/[ll]|{i^{!}}\mathcal{B}
\ar[rr]|{j^{\ast}}&&\ar@/_1pc/[ll]|{j_{!}}\ar@/^1pc/[ll]|{j_{\ast}}\mathcal{C}}$$
 If $(\mathcal{T}_{1},\mathcal{F}_{1})$ and $(\mathcal{T}_{2},\mathcal{F}_{2})$ are torsion pairs in $\mathcal{A}$ and $\mathcal{C}$, respectively. Define
$$\mathcal{T}=\{B\in \mathcal{B}~|~i^{\ast }B\in\mathcal{T}_{1}~\text{and}~j^{\ast}B\in \mathcal{T}_{2}  \}
~~\mbox{and}~~
\mathcal{F}=\{B\in \mathcal{B}~|~i^{!}B\in\mathcal{F}_{1}~\text{and}~j^{\ast}B\in \mathcal{F}_{2}  \}.$$
Then the following statements hold.
\begin{itemize}

\item[\rm (1)] \label{3} If $i^{!}$, $i{^{\ast}}$ are exact, then $(\mathcal{T},\mathcal{F})$ is a torsion pair in $\mathcal{B}$.
 We call $(\mathcal{T},\mathcal{F})$  is ``glued" by $(\mathcal T_1,\mathcal F_1)$
and $(\mathcal T_2,\mathcal F_2)$;

\item[\rm (2)] $(\mathcal{T}_{1}, \mathcal{F}_{1})=(i^{\ast }\mathcal{T}, i^{!}\mathcal{F})$ and $(\mathcal{T}_{2}, \mathcal{F}_{2})=(j^{\ast }\mathcal{T}, j^{\ast }\mathcal{F}).$
\end{itemize}
\begin{proof}
 $(1)$ Let $X\in \mathcal{T}$ and $Y\in \mathcal{F}$. Since $i^{\ast}$ is exact, there is an $\mathbb{E}_\mathcal{B}$-triangle by Lemma \ref{1} $$\xymatrix{ j_! j^\ast X\ar[r]^-{}&X\ar[r]^-{}&i_\ast i^\ast X \ar@{-->}[r]&}.$$ Applying the functor $\Hom_{\mathcal{B}}(-, Y)$ to the above $\mathbb{E}_\mathcal{B}$-triangle, we get an exact sequence
 $$\xymatrix{ \Hom_{\mathcal{B}}(i_\ast i^\ast X, Y)\ar[r]^-{}&\Hom_{\mathcal{B}}(X, Y)\ar[r]^-{}&\Hom_{\mathcal{B}}(j_! j^\ast X, Y)}.$$
 By assumption, $(\mathcal{T}_{1},\mathcal{F}_{1})$ and $(\mathcal{T}_{2},\mathcal{F}_{2})$ are torsion pairs in $\mathcal{A}$ and $\mathcal{C}$, respectively. Since $i^{\ast }X\in\mathcal{T}_{1},~ i^{!}Y\in\mathcal{F}_{1},~ j^{\ast }X\in\mathcal{T}_{2},$ and $j^{\ast }Y\in\mathcal{F}_{2}$, we have
 $$\Hom_{\mathcal{B}}(j_! j^\ast X, Y)\cong\Hom_{\mathcal{C}}( j^\ast X, j^\ast Y)=0,$$
 $$\Hom_{\mathcal{B}}(i_\ast i^\ast X, Y)\cong\Hom_{\mathcal{A}}( i^\ast X, i_! Y)=0.$$
 It follows that $\Hom_{\mathcal{B}}(X, Y)=0,$ therefore, $ \Hom_{\mathcal{B}}( \mathcal{T},\mathcal{F})=0.$

For any $B\in \mathcal{B}$, there exists an $\mathbb{E}_\mathcal{C}$-triangle $T_{2}\stackrel{}{\longrightarrow}j^{\ast}B\stackrel{}{\longrightarrow}F_{2}\stackrel{}\dashrightarrow$ with $F_{2}\in \mathcal{F}_{2}$ and $T_{2}\in \mathcal{T}_{2}$, since $j^{\ast}B\in \mathcal{C}$ and $(\mathcal{T}_{2},\mathcal{F}_{2})$ is a torsion pair in $\mathcal{C}$. Since $i^{!}$ is exact, by Lemma \ref {CY} $(8')$, $j_{\ast}$ is exact. Applying $j_{\ast}$ to the above $\mathbb{E}_\mathcal{C}$-triangle, we obtain an $\mathbb{E}_\mathcal{B}$-triangle $$j_{\ast}T_{2}\stackrel{}{\longrightarrow}j_{\ast}j^{\ast}B\stackrel{}{\longrightarrow}j_{\ast}F_{2}\stackrel{}\dashrightarrow.$$ Notice that $i^{!}$ is exact and $B\in \mathcal{B}$, there is an $\mathbb{E}_\mathcal{B}$-triangle by Lemma \ref{1} $$i_\ast i^! B\stackrel{}{\longrightarrow}B\stackrel{}{\longrightarrow}j_{\ast}j^{\ast}B\stackrel{}\dashrightarrow.$$ By $\rm (ET4)^{op}$, we have the following exact commutative diagram
\begin{equation}
\label{7}
\begin{array}{l}
\xymatrix{
  i_{\ast}i^!B \ar@{=}[d] \ar[r] & E  \ar[d] \ar[r] & j_{\ast}T_{2}  \ar[d] \\
   i_{\ast}i^!B \ar[r] &B \ar[d] \ar[r] & j_{\ast}j^{\ast}B \ar[d] \\
  & j_{\ast}F_{2} \ar@{=}[r] &j_{\ast}F_{2}. }
   \end{array}
\end{equation}
Since $i^{\ast}E\in \mathcal{A}$ and $(\mathcal{T}_{1},\mathcal{F}_{1})$ is a torsion pair in $\mathcal{A}$, there exists an $\mathbb{E}_\mathcal{A}$-triangle $$T_{1}\stackrel{}{\longrightarrow}i^{\ast}E\stackrel{}{\longrightarrow}F_{1}\stackrel{}\dashrightarrow$$ with $F_{1}\in \mathcal{F}_{1}$ and $T_{1}\in \mathcal{T}_{1}$. Because $i_{\ast}$ is exact, applying $i_{\ast}$ to the above $\mathbb{E}_\mathcal{A}$-triangle, we obtain an $\mathbb{E}_\mathcal{B}$-triangle $i_{\ast}T_{1}\stackrel{}{\longrightarrow}i_{\ast}i^{\ast}E\stackrel{}{\longrightarrow}i_{\ast}F_{1}\stackrel{}\dashrightarrow$. Notice that $i^{\ast}$ is exact, there is an $\mathbb{E}_\mathcal{B}$-triangle by Lemma \ref{1} $ j_{!}j_\ast E\stackrel{}{\longrightarrow}E\stackrel{}{\longrightarrow}i_{\ast}i^{\ast}E\stackrel{}\dashrightarrow$. By $\rm (ET4)^{op}$, we have the following exact commutative diagram
\begin{equation}\label{7}
\begin{array}{l}
\xymatrix{
   j_{!}j_\ast E \ar@{=}[d] \ar[r] &T \ar[d] \ar[r] & i_{\ast}T_{1}  \ar[d] \\
   j_{!}j_\ast E  \ar[r] & E \ar[d] \ar[r] & i_{\ast}i^{\ast}E\ar[d] \\
  &  i_{\ast}F_{1} \ar@{=}[r] & i_{\ast}F_{1}. }
  \end{array}
\end{equation}
Considering the second column of (3.1) and the second column of (3.2). By $\rm (ET4)$, we have an exact commutative diagram
\begin{equation}\label{7}
\begin{array}{l}
\xymatrix{
  T \ar@{=}[d] \ar[r] & E \ar[d] \ar[r] & i_{\ast}F_{1}  \ar[d] \\
 T \ar[r] & B\ar[d] \ar[r] & F \ar[d] \\
  & j_\ast F_{2}  \ar@{=}[r] &j_\ast F_{2}. }
   \end{array}
\end{equation}
Thus we obtain an $\mathbb{E}_\mathcal{B}$-triangle by (3.3) $$\xymatrix{ T\ar[r]^-{}&B\ar[r]^-{}& F\ar@{-->}[r]&}.$$
To get the assertion, it suffices to show $T\in \mathcal{T}$ and $F\in \mathcal{F}$.

Since $i^{\ast}$ is exact, $i^{\ast}j_{!}=0$ and $i^{\ast}i_{\ast}\Rightarrow\rm{\Id}_\mathcal{A}$ is a natural isomorphism, applying the functor $i^{\ast}$ to the first row in diagram (3.2), we have
$$i^{\ast}T\cong i^{\ast}i_{\ast}T_{1}\cong T_{1} \in \mathcal{T}_{1}$$
by Lemma \ref{00}. Since $j^{\ast}$ is exact, $j^{\ast}i_{\ast}=0$ and $j^{\ast}j_{\ast}\Rightarrow\rm{\Id}_\mathcal{C}$ is a natural isomorphism, applying the functor $j^{\ast}$ to the first row in diagram (3.1), we have $$j^{\ast}E\cong j^{\ast}j_{\ast}T_{2}\cong T_{2} \in \mathcal{T}_{2}$$
by Lemma \ref{00}. On the other hand, applying the functor $j^{\ast}$ to the first row in diagram (3.3),  we have $$j^{\ast}T\cong j^{\ast}E\in \mathcal{T}_{2}$$
by Lemma \ref{00},
It implies $T\in \mathcal{T}.$

Since $i^{!}$ is exact, $i^{!}j_{\ast}=0$ and ${\rm{\Id}}_\mathcal{A}\Rightarrow i^{!}i_{\ast}$ is a natural isomorphism, applying the functor $i^{!}$ to the third column in diagram (3.3), we have
$$i^{!}F\cong i^{!}i_{\ast}F_{1}\cong F_{1} \in \mathcal{F}_{1}$$
by Lemma \ref{00}.  Applying the functor $j^{\ast}$ to the third column in diagram (3.3), we have
$$j^{\ast}F\cong j^{\ast}j_{\ast}F_{2}\cong F_{2} \in \mathcal{F}_{2}$$
by Lemma \ref{00}. It implies $F\in \mathcal{F}.$

This show that $(\mathcal{T},\mathcal{F})$ is a torsion pair in $\mathcal{B}$.

$(2)$ It is obvious that $i^{\ast}\mathcal{T}\subseteq\mathcal{T}_{1}.$ On the other hand, for any ${T}_{1}\in\mathcal{T}_{1},$ since $$i^{\ast}i_{\ast}{T}_{1}\cong{T}_{1}\in\mathcal{T}_{1},~~~ j^{\ast}i_{\ast}{T}_{1}=0\in\mathcal{T}_{2},$$
we have $\mathcal {T}_{1}\subseteq i^{\ast}\mathcal{T}.$
Similarly, we obtain $\mathcal{F}_{1}=i^{!}\mathcal{F},$   $\mathcal{T}_{2}=j^{\ast }\mathcal{T},  \mathcal{F}_{2}=j^{\ast }\mathcal{F}.$

 This shows that $(\mathcal{T}_{1}, \mathcal{F}_{1})=(i^{\ast }\mathcal{T}, i^{!}\mathcal{F})$ and $(\mathcal{T}_{2}, \mathcal{F}_{2})=(j^{\ast }\mathcal{T}, j^{\ast }\mathcal{F}).$
\end{proof}
\end{theorem}
By applying Theorem \ref{main} to triangulated categories, we have the following.

\begin{corollary}\rm Assume that ($\mathcal{A}$, $\mathcal{B}$, $\mathcal{C}$) is a recollement of triangulated categories. Let $(\mathcal{T}_{1},\mathcal{F}_{1})$ and $(\mathcal{T}_{2},\mathcal{F}_{2})$ be torsion pairs in $\mathcal{A}$ and $\mathcal{C}$, respectively. Set
$\mathcal{T}=\{B\in \mathcal{B}~|~i^{\ast }B\in\mathcal{T}_{1}~\text{and}~j^{\ast}B\in \mathcal{T}_{2}  \}$
and $\mathcal{F}=\{B\in \mathcal{B}~|~i^{!}B\in\mathcal{F}_{1}~\text{and}~j^{\ast}B\in \mathcal{F}_{2}  \}.$ Then

$(1)$ $(\mathcal{T},\mathcal{F})$ is a torsion pair in $\mathcal{B}$;

$(2)$ $(\mathcal{T}_{1}, \mathcal{F}_{1})=(i^{\ast }\mathcal{T}, i^{!}\mathcal{F}),$ and $(\mathcal{T}_{2}, \mathcal{F}_{2})=(j^{\ast }\mathcal{T}, j^{\ast }\mathcal{F}).$
\begin{proof}This follows from Theorem \ref{main}, Remark \ref{111}, Remark \ref{77} and Remark \ref{222}.
\end{proof}
\end{corollary}

By applying Theorem \ref{main} to abelian categories, we have the following.

\begin{corollary}\rm{\cite[Theorem 1]{MH}} Assume that ($\mathcal{A}$, $\mathcal{B}$, $\mathcal{C}$) is a recollement of abelian categories. Let $(\mathcal{T}_{1},\mathcal{F}_{1})$ and $(\mathcal{T}_{2},\mathcal{F}_{2})$ be torsion pairs in $\mathcal{A}$ and $\mathcal{C}$, respectively. Set
$\mathcal{T}=\{B\in \mathcal{B}~|~i^{\ast }B\in\mathcal{T}_{1}~\text{and}~j^{\ast}B\in \mathcal{T}_{2}  \}$
and $\mathcal{F}=\{B\in \mathcal{B}~|~i^{!}B\in\mathcal{F}_{1}~\text{and}~j^{\ast}B\in \mathcal{F}_{2}  \}.$ Then:

$(1)$ $(\mathcal{T},\mathcal{F})$ is a torsion pair in $\mathcal{B}$;

$(2)$ $(\mathcal{T}_{1}, \mathcal{F}_{1})=(i^{\ast }\mathcal{T}, i^{!}\mathcal{F}),$ and $(\mathcal{T}_{2}, \mathcal{F}_{2})=(j^{\ast }\mathcal{T}, j^{\ast }\mathcal{F}).$
\begin{proof} We know that this assumption of  $i^{!}$, $i{^{\ast}}$ are exact is not necessary in abelian categories from the proof of Theorem 3.4 and Lemma 1(4) in \cite{MH}. This follows from Theorem \ref{main}, Remark \ref{111}, Remark \ref{77} and Remark \ref{222}.
\end{proof}
\end{corollary}

Our second main result shows that the converse of Theorem \ref{main} (1) holds true under certain conditions.

\begin{theorem}\label{main1}\rm
Let $\mathcal A,\mathcal B,\mathcal C$ be three extriangulated category. Assume that
 $\mathcal B$ admits a recollement relative to $\mathcal A$ and $\mathcal C$ as follows
$$\xymatrix{\mathcal{A}\ar[rr]|{i_{*}}&&\ar@/_1pc/[ll]|{i^{*}}\ar@/^1pc/[ll]|{i^{!}}\mathcal{B}
\ar[rr]|{j^{\ast}}&&\ar@/_1pc/[ll]|{j_{!}}\ar@/^1pc/[ll]|{j_{\ast}}\mathcal{C}}$$
If $(\mathcal{T},\mathcal{F})$ is a torsion pairs in $\mathcal{B}$, then
the following statements hold.

$(1)$ If $i^{!}$ is exact, $i_{\ast}i^{!}\mathcal{T}\subseteq \mathcal{T}$ and $i_{\ast}i^{\ast}\mathcal{T}\subseteq\mathcal{T}$, then $(i^{\ast}\mathcal{T}, i^{!}\mathcal{F})$ is a torsion pair in $\mathcal{A}$.

$(2)$ If $j_{!}j^{\ast}\mathcal{T}\subseteq \mathcal{T}$ or $j_{\ast}j^{\ast}\mathcal{F}\subseteq\mathcal{F}$, then $(j^{\ast}\mathcal{T}, j^{\ast}\mathcal{F})$ is a torsion pair in $\mathcal{C}$.

\begin{proof}
$(1)$ For any $A\in\mathcal{A}$, $i^{\ast}A\in\mathcal{B}$, there exists an $\mathbb{E}_\mathcal{B}$-triangle $$\xymatrix{ T\ar[r]^-{}&i^{\ast}A\ar[r]^-{}& F\ar@{-->}[r]&}$$ with $T\in \mathcal{T}$ and $F\in\mathcal{F}$ since $(\mathcal{T},\mathcal{F})$ be a torsion pairs in $\mathcal{B}$. Because $i^{!}$ is exact and $\Id_{\mathcal{A}}\Rightarrow i^{!}i_{\ast}$ is a natural isomorphism, applying $i^{!}$ to the above $\mathbb{E}_\mathcal{B}$-triangle, we obtain an $\mathbb{E}_\mathcal{A}$-triangle $$\xymatrix{  i^{!}T\ar[r]^-{}&A\ar[r]^-{}& i^{!}F\ar@{-->}[r]&}.$$
Since $i_{\ast}i^{!}{T}\in i_{\ast}i^{!}\mathcal{T}\subseteq\mathcal{T}$, note that $i^{\ast}i_{\ast}\Rightarrow\Id_{\mathcal{A}}$ is a natural isomorphism, we obtain that $i^{!}{T}\in i^{\ast}\mathcal{T}$.

Moreover, for any $T\in\mathcal{T}$ and $F\in\mathcal{F}$, we have $$\xymatrix{ \Hom_{\mathcal{A}}(i^{\ast}T, i^{!}F)\cong\Hom_{\mathcal{B}}(i_{\ast}i^{\ast}T, F)\cong\Hom_{\mathcal{B}}(T^{\prime}, F)}=0,$$ for some $T^{\prime}\in\mathcal{T}$. This shows that $(i^{\ast}\mathcal{T}, i^{!}\mathcal{F})$ is a torsion pair in $\mathcal{A}$.

$(2)$ The proof is similar to $(1)$.
\end{proof}

\end{theorem}
By applying Theorem \ref{main1} to triangulated categories, we have the following.

\begin{corollary}\label{main2}{\rm \cite[Theorem 3.3]{C}}
\rm Assume that ($\mathcal{A}$, $\mathcal{B}$, $\mathcal{C}$) is a recollement of triangulated categories. Let $(\mathcal{T},\mathcal{F})$ be a torsion pairs in $\mathcal{B}$. Then

$(1)$ If $i_{\ast}i^{!}\mathcal{T}\subseteq \mathcal{T}$ and $i_{\ast}i^{\ast}\mathcal{T}\subseteq\mathcal{T}$, then $(i^{\ast}\mathcal{T}, i^{!}\mathcal{F})$ is a torsion pair in $\mathcal{A}$.

$(2)$ If $j_{!}j^{\ast}\mathcal{T}\subseteq \mathcal{T}$ or $j_{\ast}j^{\ast}\mathcal{F}\subseteq\mathcal{F}$, then $(j^{\ast}\mathcal{T}, j^{\ast}\mathcal{F})$ is a torsion pair in $\mathcal{C}$.

\begin{proof}
This follows from Theorem \ref{main1}, Remark \ref{111}, Remark \ref{77} and Remark \ref{222}.
\end{proof}
\end{corollary}

By applying Theorem \ref{main1} to abelian categories, we have the following.

\begin{corollary}\label{main22}\rm{\cite[Theorem 2]{MH}} Assume that ($\mathcal{A}$, $\mathcal{B}$, $\mathcal{C}$) is a recollement of abelian categories. Let $(\mathcal{T},\mathcal{F})$ be a torsion pairs in $\mathcal{B}$. Then

$(1)$ $(i^{\ast}\mathcal{T}, i^{!}\mathcal{F})$ is a torsion pair in $\mathcal{A}$.

$(2)$ If $j_{!}j^{\ast}\mathcal{T}\subseteq \mathcal{T}$ or $j_{\ast}j^{\ast}\mathcal{F}\subseteq\mathcal{F}$, then $(j^{\ast}\mathcal{T}, j^{\ast}\mathcal{F})$ is a torsion pair in $\mathcal{C}$.

\begin{proof}
This follows from Lemma 1, 3 in \cite{MH}, Theorem \ref{main1}, Remark \ref{111}, Remark \ref{77} and Remark \ref{222}.
\end{proof}
\end{corollary}

\section{From recollement of extriangulated categories to recollement of abelian categories }
In this section, when we say that $\mathcal{T}$ is a subcategory of $\mathcal{C}$, we always mean that $\mathcal{T}$ is a full subcategory which is closed under isomorphisms.
\begin{definition}\label{dd1}\rm{\cite[Definition 3.21]{ZZ1}}
Let $\mathcal{C}$ be an  extriangulated category. A subcategory $\mathcal{T}$ of $\mathcal{C}$ is called
\emph{strongly contravariantly finite}, if for any object $C\in\mathcal{C}$, there exists an $\E$-triangle
$$\xymatrix{K\ar[r]&T\ar[r]^{g}&C\ar@{-->}[r]^{\del}&,}$$
where $g$ is a right $\mathcal{T}$-approximation of $C$.

Dually, a subcategory $\mathcal{T}$ of $\mathcal{C}$ is called
\emph{strongly  covariantly  finite}, if for any object $C\in\mathcal{C}$, there exists an $\E$-triangle
$$\xymatrix{C\ar[r]^{f}&T\ar[r]&L\ar@{-->}[r]^{\del'}&,}$$
where $f$ is a left $\mathcal{T}$-approximation of $C$.

A strongly contravariantly finite and strongly  covariantly finite subcategory is called \emph{ strongly functorially finite}.
\end{definition}

\begin{definition}\label{y5}\rm{\cite[Definition 3.1]{CZZ}}
Let $\mathcal{C}$ be an extriangulated category, $\mathcal{T}$ a subcategory of $\mathcal{C}$.
\begin{itemize}

\item $\mathcal{T}$  is called rigid if $\E(\mathcal{T}, \mathcal{T})=0$.

\item $\mathcal{T}$ is called \emph{cluster tilting} if it satisfies the following conditions:
\begin{enumerate}
\item[(1)] $\mathcal{T}$ is strongly functorially finite in $\mathcal{C}$;
\item[(2)] $T\in \mathcal{T}$ if and only if $\E(T, \mathcal{T})=0$;
\item[(3)] $T\in \mathcal{T}$ if and only if $\E(\mathcal{T}, T)=0$.
\end{enumerate}

\end{itemize}
\end{definition}
By definition of a cluster tilting subcategory, we can immediately conclude as follows.

\begin{remark}\label{y6}
Let $\mathcal{C}$ be an extriangulated category and $\mathcal{T}$ a subcategory of $\mathcal{C}$.
\begin{itemize}
\item $\mathcal{T}$ is a cluster tilting subcategory of $\mathcal{C}$, if and only if
\begin{enumerate}
\item[(1)] $\mathcal{T}$ is rigid;
\item[(2)] For any object $C\in\mathcal{C}$, there exists an $\E$-triangle
$\xymatrix{C\ar[r]&T_{1}\ar[r]&T_{2}\ar@{-->}[r]^{\del'}&,}$
where $T_{1}, T_{2}\in\mathcal{T}$.
\item[(3)]  For any object $C\in\mathcal{C}$, there exists an $\E$-triangle
$\xymatrix{T_{3}\ar[r]&T_{4}\ar[r]&C\ar@{-->}[r]^{\eta'}&,}$
where $T_{3}, T_{4}\in\mathcal{T}$.
\end{enumerate}

\end{itemize}
\end{remark}

\begin{lemma}\label{hz}\rm{\cite[Theorem 3.4]{ZZ2} and \cite[Theorem 3.3]{HZ}}
Let $\mathcal{C}$ be a extriangulated categories and $\mathcal{T}$ a cluster-tilting subcategory of $\mathcal{C}$. Then $\mathcal{C}/\mathcal{T}$ is an abelian category.

\end{lemma}

Our third main result is the following.

\begin{theorem}\label{zh}
Let $\mathcal{A}$, $\mathcal{B}$ and $\mathcal{C}$ be three extriangulated categories, where $\mathcal{B}$ has enough projectives. Assume that $\mathcal{B}$ admits a recollement relative to $\mathcal{A}$ and $\mathcal{C}$, i.e.
\begin{equation}\label{recol}
  \xymatrix{\mathcal{A}\ar[rr]|{i_{*}}&&\ar@/_1pc/[ll]|{i^{*}}\ar@/^1pc/[ll]|{i^{!}}\mathcal{B}
\ar[rr]|{j^{\ast}}&&\ar@/_1pc/[ll]|{j_{!}}\ar@/^1pc/[ll]|{j_{\ast}}\mathcal{C}}.
\end{equation}
If ${i}^{*}$ and $ {i}^{!}$ are exact, $\mathcal{T}$ is a \emph{cluster tilting} subcategory of $\mathcal{B}$ and satisfies ${j_{\ast}j^{\ast}}\mathcal{T}\subseteq\mathcal{T}, {i_{\ast}i^{\ast}}\mathcal{T}\subseteq\mathcal{T}$. Then abelian category $\mathcal{B}/\mathcal{T}$ admits a recollement relative to abelian category $\mathcal{A}/i^{\ast}\mathcal{T}$ and $\mathcal{C}/j^{\ast}\mathcal{T}$ as follows:
\begin{equation}\label{recoll}
  \xymatrix{\mathcal{A}/i^{\ast}\mathcal{T}\ar[rr]|{\overline{i}_{*}}&&\ar@/_1pc/[ll]|{\overline{i}^{*}}\ar@/^1pc/[ll]|{\overline{i}^{!}}\mathcal{B}/\mathcal{T}
\ar[rr]|{\overline{j}^{\ast}}&&\ar@/_1pc/[ll]|{\overline{j}_{!}}\ar@/^1pc/[ll]|{\overline{j}_{\ast}}\mathcal{C}/j^{\ast}\mathcal{T}}.
\end{equation}

\end{theorem}

In order to prove Theorem \ref{zh}, we need some preparations as follows.

\begin{lemma}\label{OO}
Let $(\mathcal{A}, \mathcal{B},\mathcal{C})$ be a recollement of extriangulated categories. Assume that $\mathcal{B}$ has enough projectives, ${i}^{*}$ and $ {i}^{!}$ are exact. If $\mathcal{T}$ is a cluster tilting subcategory of $\mathcal{B}$ and satisfies ${j_{\ast}j^{\ast}}\mathcal{T}\subseteq\mathcal{T}$, ${i_{\ast}i^{\ast}}\mathcal{T}\subseteq\mathcal{T}$. Then

$(1)$ ${j^{\ast}}\mathcal{T}$ is a cluster tilting subcategory of $\mathcal{C}$;

$(2)$ ${i^{\ast}}\mathcal{T}$ is a cluster tilting subcategory of $\mathcal{A}$.
\begin{proof} $(1)$ For any $C\in\mathcal{C}$, $j_{\ast}C\in\mathcal{B}$, since $\mathcal{T}$ is a cluster tilting subcategory of $\mathcal{B}$, then there exists an $\mathbb{E}_\mathcal{B}$-triangle $$\xymatrix{ j_{\ast}C\ar[r]^-{}&T_1\ar[r]^-{}& T_2\ar@{-->}[r]&}$$ with $T_1, T_2\in \mathcal{T}$. Because $j^{\ast}$ is exact and $j^{\ast}j_{\ast}\Rightarrow \Id_{\mathcal{C}}$ is a natural isomorphism, applying $j^{\ast}$ to the above $\mathbb{E}_\mathcal{B}$-triangle, we obtain an $\mathbb{E}_\mathcal{C}$-triangle $$\xymatrix{ C\ar[r]^-{}&j^{\ast}T_1\ar[r]^-{}& j^{\ast}T_2\ar@{-->}[r]&}$$ with $j^{\ast}T_1, j^{\ast}T_2\in j^{\ast}\mathcal{T}$.

Similarly, for any $C\in\mathcal{C}$, we have an $\mathbb{E}_\mathcal{C}$-triangle $$\xymatrix{ j^{\ast}T_3\ar[r]^-{}&j^{\ast}T_4\ar[r]^-{}&C \ar@{-->}[r]&}$$ with $j^{\ast}T_3, j^{\ast}T_4\in j^{\ast}\mathcal{T}$.

For any $X_1, X_2 \in\mathcal{T}$, since ${i}^{!}$ is exact,  ${j}^{*}$ preserves projective objects by Lemma \ref{CY} (4), $(8')$. Since $\mathcal{B}$ has enough projectives, we have $\mathbb{E}_{\mathcal{C}}(j^{\ast}X_1, j^{\ast}X_2)\cong\mathbb{E}_{\mathcal{B}}(X_1, j_{\ast}j^{\ast}X_2)=0$ by Lemma 2.16 in \cite{WWZ}. i.e. $\mathbb{E}_{\mathcal{C}}(j^{\ast}\mathcal{T}, j^{\ast}\mathcal{T})=0$. So ${j^{\ast}}\mathcal{T}$ is a cluster tilting subcategory of $\mathcal{C}$ by Remark \ref{y6}.

$(2)$ It is similar to $(1)$.
\end{proof}
\end{lemma}

The following statement is well known. Let $\mathcal{C}$ is an additive category and $\mathcal{T}$ is a subcategory. For the quotient category $\mathcal{C}/\mathcal{T}$, let ${\pi: C \rightarrow \mathcal{C}/\mathcal{T}}$ be the projective functor, then ${\pi}$ has the universai property. That is, for additive category $\mathcal{B}$ and ${F: C \rightarrow \mathcal{B}}$ with ${F(T)}$ is null for any $T\in\mathcal{T}$, there exists a unique additive functor ${G: \mathcal{C}/\mathcal{T} \rightarrow \mathcal{B}}$, such that the following diagram commutes:
$$\xymatrix@R=1cm@C=1.2cm{
\mathcal{C} \ar[dr]^F \ar[d]^\pi \\
\mathcal{C}/\mathcal{T}\ar@{-->}[r]^{G} & B}$$

\begin{lemma}\label{j}
Let $(\mathcal{A}, \mathcal{B},\mathcal{C})$ be a recollement of extriangulated categories. Assume that $\mathcal{B}$ has enough projectives, ${i}^{*}$ and $ {i}^{!}$ are exact. If $\mathcal{T}$ is a cluster tilting subcategory of $\mathcal{B}$ and satisfies ${j_{\ast}j^{\ast}}\mathcal{T}\subseteq\mathcal{T}, {i_{\ast}i^{\ast}}\mathcal{T}\subseteq\mathcal{T}$. Then

$(1)$ ${j^{\ast}}$ induces an additive functor of abelian category ${\overline{j}^{\ast}:\mathcal{B}/\mathcal{T}\rightarrow \mathcal{C}/j^{\ast}\mathcal{T}}$;

$(2)$ ${j_{\ast}}$ induces an additive functor of abelian category ${\overline{j}_{\ast}: \mathcal{C}/j^{\ast}\mathcal{T}\rightarrow\mathcal{B}/\mathcal{T}}$;

$(3)$ $({\overline{j}^{\ast},\overline{j}_{\ast}})$ is an adjoint pair.

$(1')$ ${i^{\ast}}$ induces an additive functor of abelian category ${\overline{i}^{\ast}:\mathcal{B}/\mathcal{T}\rightarrow \mathcal{A}/i^{\ast}\mathcal{T}}$;

$(2')$ ${i_{\ast}}$ induces an additive functor of abelian category ${\overline{i}_{\ast}:\mathcal{A}/i^{\ast}\mathcal{T} \rightarrow\mathcal{B}/\mathcal{T}}$;

$(3')$ $({\overline{i}^{\ast},\overline{i}_{\ast}})$ is an adjoint pair.

\begin{proof}
$(1)$ Since ${j^{\ast}}\mathcal{T}$ is a cluster tilting subcategory of $\mathcal{C}$ by Lemma \ref{OO}, $\mathcal{B}/\mathcal{T}, \mathcal{C}/j^{\ast}\mathcal{T}$ are abelian categories by Lemma \ref{hz}. Since $\pi_\mathcal{C}j^{\ast}\mathcal{T}=0$ for any ${T}\in\mathcal{T}$, there exists an unique additive functor ${\overline{j}^{\ast}:\mathcal{B}/\mathcal{T}\rightarrow \mathcal{C}/j^{\ast}\mathcal{T}}$, such that the following diagram is commutative:
$$\xymatrix{
\mathcal{B} \ar[r]^{j^\ast} \ar[d]^{\pi_\mathcal{B} }&\mathcal{C}\ar[d]^{\pi_\mathcal{C}} \\
\mathcal{B}/\mathcal{T}\ar@{-->}[r]^{\overline{j}^{\ast}} &\mathcal{C}/j^{\ast}\mathcal{T}}$$

$(2)$ Since $j_{\ast}j^{\ast}\mathcal{T}\subseteq \mathcal{T}$, $\pi_\mathcal{B}j_{\ast}j^{\ast}\mathcal{T}=0$ for any ${T}\in\mathcal{T}$, there exists an unique additive functor ${\overline{j}_{\ast}: \mathcal{C}/j^{\ast}\mathcal{T}\rightarrow\mathcal{B}/\mathcal{T}}$, such that the following diagram is commutative:
$$\xymatrix{
\mathcal{C}\ar[r]^{j_\ast} \ar[d]^{\pi_\mathcal{C}}&\mathcal{B}\ar[d]^{\pi_{\mathcal{B}}} \\
\mathcal{C}/j^{\ast}\mathcal{T}\ar@{-->}[r]^{\overline{j}_{\ast}} &\mathcal{B}/\mathcal{T}}$$

$(3)$ It is similar to the proof of Lemma 3.4 in \cite{LW}.
\end{proof}
\end{lemma}

{\bf Now we are ready to prove Theorem \ref{zh}.}

\proof By Lemma \ref{OO} and Lemma \ref{hz}, we have that $\mathcal{B}/\mathcal{T}, \mathcal{C}/j^{\ast}\mathcal{T}$ and $\mathcal{A}/i^{\ast}\mathcal{T}$ are abelian categories. By Lemma \ref{j}, ${j^{\ast}}, {j_{\ast}}, {i^{\ast}},{i_{\ast}}$ induce additive functors $\overline{j^{\ast}}, \overline{j_{\ast}}, \overline{i^{\ast}}, \overline{i_{\ast}}$, respectively. Moreover, $({\overline{j}^{\ast},\overline{j}_{\ast}}), ({\overline{i}^{\ast},\overline{i}_{\ast}})$ are adjoint pairs. Next we show that $i^{!}, i_{!}$ induce additive functors $\overline{i}^{!}, {\overline{j}_{!}}$, respectively.

Let $T \in\mathcal{T}$.  For any $T'\in\mathcal{T}$, since $\mathcal{T}$ is a cluster tilting subcategory of $\mathcal{B}$, and ${i_{\ast}i^{\ast}}\mathcal{T}\subseteq\mathcal{T}$,  we have $$\mathbb{E}_{\mathcal{B}}(T', i_{\ast}i^{!}T)\cong\mathbb{E}_{\mathcal{A}}(i^{\ast}T', i^{!}T)\cong\mathbb{E}_{\mathcal{B}}(i_{\ast}i^{\ast}T', T)=0$$ by Lemma \ref{CY} (4), $(7)$ and Lemma 2.16 in \cite{WWZ}.
Thus $i_{\ast}i^{!}T\in\mathcal{T}$. Note that $i^{!}T\cong i^{\ast}i_{\ast}i^{!}T\in i^{\ast}\mathcal{T}$, so $i^{!}$ induce additive functors $\overline{{i}^{!}}$.

Let $T \in\mathcal{T}$.  For any $T'\in\mathcal{T}$, since $j^{\ast}\mathcal{T}$ is a cluster tilting subcategory of $\mathcal{C}$, then we have $$\mathbb{E}_{\mathcal{B}}(j_{!}j^{\ast}T, T')\cong\mathbb{E}_{\mathcal{C}}(j^{\ast}T, j^{\ast}T')=0.$$
by Lemma \ref{CY} $(7')$. Thus $j_{!}j^{\ast}T\in\mathcal{T}$ which implies that $j_{!}$ induces additive functors $\overline{ j_{!}}$.

Similar to the proof of Lemma \ref{j} (3), we obtain $({\overline{i}_{\ast},\overline{i}^{!}}), ({\overline{j}_{!},\overline{j}^{\ast}})$ are adjoint pair.

Since $\Id_{\mathcal A}\Rightarrow i^{!}i_{\ast}$ is an natural isomorphism, then $\Id_{\mathcal{A}/i^{\ast}\mathcal{T} }\Rightarrow \overline i^{!}\overline i_{\ast}$ is an natural isomorphism,. It follows that $\overline i_{\ast}$ is fully faithful. Similarly, we get that $\overline j_{\ast}, \overline j_{!}$ are fully faithful.

It is easy to see $\overline j^{\ast}\overline i_{\ast}=0$ since $j^{\ast} i_{\ast}=0$.

This completes the proof.  \qed
\medskip

By applying Theorem \ref{zh} to triangulated categories, we get the following.

\begin{corollary}\rm{\cite[Theorem 1.1]{LW}}
Let $(\mathcal{A}, \mathcal{B}, \mathcal{C})$ be a recollement of triangulated categories. If $\mathcal{T}$ is a cluster tilting subcategory of $\mathcal{B}$ with ${j_{\ast}j^{\ast}}\mathcal{T}\subseteq\mathcal{T}, {i_{\ast}i^{\ast}}\mathcal{T}\subseteq\mathcal{T}$, then the abelian category $\mathcal{B}/\mathcal{T}$ admits a recollement relative to abelian category $\mathcal{A}/i^{\ast}\mathcal{T}$ and $\mathcal{C}/j^{\ast}\mathcal{T}$ as follows:

 $$ \xymatrix{\mathcal{A}/i^{\ast}\mathcal{T}\ar[rr]|{\overline{i}_{*}}&&\ar@/_1pc/[ll]|{\overline{i}^{*}}\ar@/^1pc/[ll]|{\overline{i}^{!}}\mathcal{B}/\mathcal{T}
\ar[rr]|{\overline{j}^{\ast}}&&\ar@/_1pc/[ll]|{\overline{j}_{!}}\ar@/^1pc/[ll]|{\overline{j}_{\ast}}\mathcal{C}/j^{\ast}\mathcal{T}}.$$
\begin{proof} Since a triangulated category always has enough projectives. This follows from Theorem \ref{zh} and Remark \ref{111}.
\end{proof}
\end{corollary}

\section{An example}
In this section, we give an example to explain our main results.

\begin{example}\label{example1}
Let $A$ be the path algebra of the quiver $1\rightarrow 2$ over a field. The
Auslander-Reiten quiver of the left $A$-modules category $\textrm{mod} A$ is as follows
$$\xymatrix@C=1em@R=1em{ &P(1)\ar[dr]&\\
S(2)\ar[ur]&&S(1).}$$
Then the upper triangular matrix algebra $\Lambda=\begin{bmatrix}\begin{smallmatrix} A & A\\0&A\end{smallmatrix}\end{bmatrix}
$ is given by the following quiver
$$
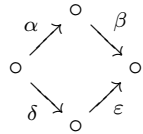
\xymatrix@C=1em@R=1em{ &\circ\ar[dr]^{\beta}&\\
\circ\ar[dr]_{\delta}\ar[ur]^{\alpha}&&\circ\\
&\circ\ar[ur]_{\varepsilon}&}$$
with the relation $\beta\alpha=\gamma\delta$. It is well-known that each left $\Lambda$-module  can be uniquely described as $\begin{bmatrix}\begin{smallmatrix} X \\Y\end{smallmatrix}\end{bmatrix}_{f}$, where $f:Y\rightarrow X$ is a left $A$-module homomorphism. Then the Auslander-Reiten quiver of $\textrm{mod} \Lambda$ is given by
 $$\xymatrix{ &\text{$\begin{bmatrix}\begin{smallmatrix} P(1) \\\textrm{O}\end{smallmatrix}\end{bmatrix}$}_{0}\ar[dr]&&\text{$\begin{bmatrix}\begin{smallmatrix} \textrm{O} \\S(2)\end{smallmatrix}\end{bmatrix}$}_{0}\ar[dr]&&\text{$\begin{bmatrix}\begin{smallmatrix} S(1) \\S(1)\end{smallmatrix}\end{bmatrix}$}_{1}\ar[dr]&\\
 \text{$\begin{bmatrix}\begin{smallmatrix} S(2) \\\textrm{O}\end{smallmatrix}\end{bmatrix}$}_{0}\ar[ur]\ar[dr]&&\text{$\begin{bmatrix}\begin{smallmatrix} P(1) \\ S(2)\end{smallmatrix}\end{bmatrix}$}_{f}\ar[r]\ar[ur]\ar[dr]&\text{$\begin{bmatrix}\begin{smallmatrix} P(1) \\ P(1)\end{smallmatrix}\end{bmatrix}$}_{1}\ar[r]&\text{$\begin{bmatrix}\begin{smallmatrix} S(1) \\ P(1)\end{smallmatrix}\end{bmatrix}$}_{g}\ar[ur]\ar[dr]&&\text{$\begin{bmatrix}\begin{smallmatrix} \textrm{O} \\S(1)\end{smallmatrix}\end{bmatrix}$}_{0}\\
 &\text{$\begin{bmatrix}\begin{smallmatrix} S(2) \\S(2)\end{smallmatrix}\end{bmatrix}$}_{0}\ar[ur]&&\text{$\begin{bmatrix}\begin{smallmatrix} S(1) \\\textrm{O}\end{smallmatrix}\end{bmatrix}$}_{0}\ar[ur]&&\text{$\begin{bmatrix}\begin{smallmatrix} \textrm{O}\\P(1) \end{smallmatrix}\end{bmatrix}$}_{0}\ar[ur]&}$$
 By \cite[Example 2.12]{P}, we have a recollement of module categories
\begin{equation}\label{eqEr1}
\xymatrix{\textrm{mod} A\ar[rr]|{i_{*}}&&\ar@/_1pc/[ll]|{i^{*}}\ar@/^1pc/[ll]|{i^{!}}\textrm{mod} \Lambda
\ar[rr]|{j^{\ast}}&&\ar@/_1pc/[ll]|{j_{!}}\ar@/^1pc/[ll]|{j_{\ast}}\textrm{mod} A.}
\end{equation}
where
\begin{align*}
i^{\ast}\big(\begin{bmatrix}\begin{smallmatrix} X \\Y\end{smallmatrix}\end{bmatrix}_{f}\big)&=\textrm{Coker} f & i_{\ast}(X)&=\begin{bmatrix}\begin{smallmatrix} X \\\textrm{O}\end{smallmatrix}\end{bmatrix}_{0}& i^{!}\big(\begin{bmatrix}\begin{smallmatrix} X \\Y\end{smallmatrix}\end{bmatrix}_{f}\big)&=X\\
j_!(Y)& =\begin{bmatrix}\begin{smallmatrix} Y \\Y\end{smallmatrix}\end{bmatrix}_{1}
& j^{\ast}\big(\begin{bmatrix}\begin{smallmatrix} X \\Y\end{smallmatrix}\end{bmatrix}_{f}\big)&=Y
& j_\ast(Y)&=\begin{bmatrix}\begin{smallmatrix}  \textrm{O} \\Y\end{smallmatrix}\end{bmatrix}_{0}.
\end{align*}
Let \begin{align*}
      \mathcal{A}&=\textrm{add}(P(1)\oplus S(2)), \\
     \mathcal{B}&=\textrm{add}\big(\begin{bmatrix}\begin{smallmatrix} P(1) \\\textrm{O}\end{smallmatrix}\end{bmatrix}_{0}\oplus\begin{bmatrix}\begin{smallmatrix} P(1) \\P(1)\end{smallmatrix}\end{bmatrix}_{1}\oplus\begin{bmatrix}\begin{smallmatrix} S(2) \\\textrm{O}\end{smallmatrix}\end{bmatrix}_{0}\oplus\begin{bmatrix}\begin{smallmatrix} \textrm{O} \\P(1)\end{smallmatrix}\end{bmatrix}_{0}\big),\\
     \mathcal{C}&=\textrm{add}P(1).
    \end{align*}

Note that $P(1)$ is a non-zero projective object and $\mathcal{A}$ is not closed under cokernel. This implies that
$\mathcal{A}$ is neither triangulated nor abelian. Clearly, $\mathcal{A}$ is an extension-closed
subcategory of $\textrm{mod} A$. Hence, $\mathcal{A}$ is an extriangulated
category which is neither abelian nor triangulated. It is easy to see that there are only two following non-split exact sequences in $\mathcal{B}$
\begin{align*}
 0\rightarrow\begin{bmatrix}\begin{smallmatrix} P(1) \\\textrm{O}\end{smallmatrix}\end{bmatrix}_{0}\rightarrow\begin{bmatrix}\begin{smallmatrix} P(1) \\P(1)\end{smallmatrix}\end{bmatrix}_{1}\rightarrow \begin{bmatrix}\begin{smallmatrix} \textrm{O} \\P(1)\end{smallmatrix}\end{bmatrix}_{0}\rightarrow0;\\[2mm]
 0\rightarrow\begin{bmatrix}\begin{smallmatrix} \textrm{O} \\P(1)\end{smallmatrix}\end{bmatrix}_{0}\rightarrow\begin{bmatrix}\begin{smallmatrix} P(1) \\P(1)\end{smallmatrix}\end{bmatrix}_{1}\rightarrow\begin{bmatrix}\begin{smallmatrix} P(1) \\\textrm{O}\end{smallmatrix}\end{bmatrix}_{0} \rightarrow0.
\end{align*}
Thus, $\mathcal{B}$ is an extension-closed
subcategory of $\textrm{mod} \Lambda$. Moreover, $\mathcal{B}$ has a projective-injective object $\begin{bmatrix}\begin{smallmatrix} P(1) \\P(1)\end{smallmatrix}\end{bmatrix}_{1}$ and the cokernel of the homomorphism $\begin{bmatrix}\begin{smallmatrix} S(2) \\\textrm{O}\end{smallmatrix}\end{bmatrix}_{0}\rightarrow \begin{bmatrix}\begin{smallmatrix} P(1) \\\textrm{O}\end{smallmatrix}\end{bmatrix}_{0}$ does not lie in $\mathcal{B}$. Therefore, $\mathcal{B}$ is just an extriangulated
category. It is easy to check that the functor $i_{\ast}$  sends $\mathcal{A}$ to $\mathcal{B}$ and  the functor $j_{!}$, $j_{\ast}$ send $\mathcal{C}$ to $\mathcal{B}$. Moreover, the functors $i^{\ast}$, $i^!$  send $\mathcal{B}$ to $\mathcal{A}$ and  the functor  $j^{\ast}$ sends $\mathcal{B}$ to $\mathcal{C}$.
Then we have a diagram of the restricted functors as follows:
\begin{equation}\label{eqEr2}
\xymatrix{\mathcal{A}\ar[rr]|{i_{*}}&&\ar@/_1pc/[ll]|{i^{*}}\ar@/^1pc/[ll]|{i^{!}}\mathcal{B}
\ar[rr]|{j^{\ast}}&&\ar@/_1pc/[ll]|{j_{!}}\ar@/^1pc/[ll]|{j_{\ast}}\mathcal{C}}
\end{equation}
Next, we check that this is a recollement of an extriangulated category which is neither abelian nor triangulated.
\begin{enumerate}
  \item[(R1)]  As $\mathcal{A}$, $\mathcal{B}$ and $\mathcal{C}$ are full subcategories, ($i^{\ast}$, $i_{\ast}$, $i^{!}$) and ($j_{!}$, $j^{\ast}$, $j_{\ast}$) are adjoint triples.
  \item[(R2)]  Im$i_{\ast}=\textrm{add}\big(\begin{bmatrix}\begin{smallmatrix} P(1) \\\textrm{O}\end{smallmatrix}\end{bmatrix}_{0}\oplus\begin{bmatrix}\begin{smallmatrix} S(2) \\\textrm{O}\end{smallmatrix}\end{bmatrix}_{0}\big)=$Ker$j^{\ast}$.
  \item[(R3)]  Since the functors $i_{\ast}$, $j_{!}$ and $j_{\ast}$ in (\ref{eqEr1}) are fully faithful, it follows that they are also fully faithful in (\ref{eqEr2}).
  \item[(R4)]  Note that $i^{!}$ is exact. For any $\begin{bmatrix}\begin{smallmatrix} X \\Y\end{smallmatrix}\end{bmatrix}_{f}\in \mathcal{B}$, by
      \cite[Proposition 2.6]{P}, there exists an
exact sequence
$$0\rightarrow i_{\ast}i^{!}(\begin{bmatrix}\begin{smallmatrix} X \\Y\end{smallmatrix}\end{bmatrix}_{f})\rightarrow \begin{bmatrix}\begin{smallmatrix} X \\Y\end{smallmatrix}\end{bmatrix}_{f}\rightarrow j_{\ast}j^{\ast}(\begin{bmatrix}\begin{smallmatrix} X \\Y\end{smallmatrix}\end{bmatrix}_{f})\rightarrow 0$$
in $\textrm{mod} \Lambda$, which also gives a left exact $\mathbb{E}$-triangle sequence in $\mathcal{B}$
$$\begin{bmatrix}\begin{smallmatrix} X \\\textrm{O}\end{smallmatrix}\end{bmatrix}_{0}\rightarrow \begin{bmatrix}\begin{smallmatrix} X \\Y\end{smallmatrix}\end{bmatrix}_{f}\rightarrow \begin{bmatrix}\begin{smallmatrix} \textrm{O} \\Y\end{smallmatrix}\end{bmatrix}_{0}\dashrightarrow.$$
\item[(R5)] For any $\begin{bmatrix}\begin{smallmatrix} X \\Y\end{smallmatrix}\end{bmatrix}_{f}\in \mathcal{B}$, by
    \cite[Proposition 2.6]{P}, there exists an exact sequence
    $$0\rightarrow i_{\ast}(X')\rightarrow j_{!}j^{\ast}(\begin{bmatrix}\begin{smallmatrix} X \\Y\end{smallmatrix}\end{bmatrix}_{f})\rightarrow \begin{bmatrix}\begin{smallmatrix} X \\Y\end{smallmatrix}\end{bmatrix}_{f}\rightarrow i_{\ast}i^{\ast}(\begin{bmatrix}\begin{smallmatrix} X \\Y\end{smallmatrix}\end{bmatrix}_{f})\rightarrow 0$$
    where $X'\in \textrm{mod} A$. One can check that $X'\in \textrm{add}P(1)\subseteq \mathcal{A}$. Hence, this exact sequence gives a right exact $\mathbb{E}$-triangle sequence
    $$i_{\ast}(X')\rightarrow \begin{bmatrix}\begin{smallmatrix} Y \\Y\end{smallmatrix}\end{bmatrix}_{1}\rightarrow \begin{bmatrix}\begin{smallmatrix} X \\Y\end{smallmatrix}\end{bmatrix}_{f}\rightarrow \begin{bmatrix}\begin{smallmatrix} \textrm{Coker} f \\\textrm{O}\end{smallmatrix}\end{bmatrix}_{0}.$$
\end{enumerate}
It is easy to check that the functors $i^{\ast}$, $j_{!}$ are right exact, $i_{\ast}$, $j^{\ast}$  exact  and $i^{!}$, $j_{\ast}$  left exact. Therefore, the diagram (\ref{eqEr2}) is  a recollement of an extriangulated category.\par
Finally, we remark that the $i^{!}$ is exact and $i^{!}$ is not left exact. In fact, on one hand, we have $\textrm{Ext}^{1}_{\Lambda}$( $\begin{bmatrix}\begin{smallmatrix} X \\Y\end{smallmatrix}\end{bmatrix}_{f}$, $\begin{bmatrix}\begin{smallmatrix} X' \\Y'\end{smallmatrix}\end{bmatrix}_{f'}$)$=0$, for all pairs ( $\begin{bmatrix}\begin{smallmatrix} X \\Y\end{smallmatrix}\end{bmatrix}_{f}$, $\begin{bmatrix}\begin{smallmatrix} X' \\Y'\end{smallmatrix}\end{bmatrix}_{f'}$) in $\mathcal{B}$ but the pairs ($\begin{bmatrix}\begin{smallmatrix} \textrm{O}\\P(1)\end{smallmatrix}\end{bmatrix}_{0}$, $\begin{bmatrix}\begin{smallmatrix} P(1) \\\textrm{O}\end{smallmatrix}\end{bmatrix}_{0}$) and ($\begin{bmatrix}\begin{smallmatrix} P(1) \\\textrm{O}\end{smallmatrix}\end{bmatrix}_{0}$,$\begin{bmatrix}\begin{smallmatrix} \textrm{O}\\P(1)\end{smallmatrix}\end{bmatrix}_{0}$). The unique non-trivial element of $\textrm{Ext}^{1}_{\Lambda}$($\begin{bmatrix}\begin{smallmatrix} \textrm{O}\\P(1)\end{smallmatrix}\end{bmatrix}_{f}$, $\begin{bmatrix}\begin{smallmatrix} P(1) \\\textrm{O}\end{smallmatrix}\end{bmatrix}_{f'}$) is $$[0\rightarrow\begin{bmatrix}\begin{smallmatrix} P(1) \\\textrm{O}\end{smallmatrix}\end{bmatrix}_{0}\rightarrow\begin{bmatrix}\begin{smallmatrix} P(1) \\P(1)\end{smallmatrix}\end{bmatrix}_{1}\rightarrow \begin{bmatrix}\begin{smallmatrix} \textrm{O} \\P(1)\end{smallmatrix}\end{bmatrix}_{0}\rightarrow0],$$ and that of $\textrm{Ext}^{1}_{\Lambda}$($\begin{bmatrix}\begin{smallmatrix} P(1) \\\textrm{O}\end{smallmatrix}\end{bmatrix}_{0}$,$\begin{bmatrix}\begin{smallmatrix} \textrm{O}\\P(1)\end{smallmatrix}\end{bmatrix}_{0}$) is $$[0\rightarrow\begin{bmatrix}\begin{smallmatrix} \textrm{O} \\P(1)\end{smallmatrix}\end{bmatrix}_{0}\rightarrow\begin{bmatrix}\begin{smallmatrix} P(1) \\P(1)\end{smallmatrix}\end{bmatrix}_{1}\rightarrow \begin{bmatrix}\begin{smallmatrix} P(1) \\\textrm{O}\end{smallmatrix}\end{bmatrix}_{0}\rightarrow0].$$ Applying $i^{!}$ to these two sequences, we can get the conflations
\begin{align*}
 P(1)\xrightarrow{1}P(1)\rightarrow 0\dashrightarrow, \\
0\rightarrow P(1)\xrightarrow{1}P(1)\dashrightarrow.
\end{align*}
 Hence, $i^{!}$ is exact. On the other hand,  the image of the $\mathbb{E}$-triangle $\begin{bmatrix}\begin{smallmatrix} P(1) \\\textrm{O}\end{smallmatrix}\end{bmatrix}_{0}\rightarrow\begin{bmatrix}\begin{smallmatrix} P(1) \\P(1)\end{smallmatrix}\end{bmatrix}_{1}\rightarrow \begin{bmatrix}\begin{smallmatrix} \textrm{O} \\P(1)\end{smallmatrix}\end{bmatrix}_{0} \dashrightarrow $ under the functor $i^{\ast}$ is the right exact sequence $P(1)\rightarrow 0\rightarrow 0$ in $\mathcal{A}$. Thus, it is not left exact.\par

Now, let $\mathcal{T}_{1}=\textrm{add}P(1)$,  $\mathcal{F}_{1}=\textrm{add}S(2)$, $\mathcal{T}_{2}=\textrm{add}P(1)$ and $\mathcal{F}_{2}=0$. Then
($\mathcal{T}_{1}$, $\mathcal{F}_{1}$) and ($\mathcal{T}_{2}$, $\mathcal{F}_{2}$) are torsion pair in $\mathcal{A}$  and $\mathcal{C}$, respectively. By the construction in Theorem \ref{main}, we get the pair
\begin{align*}
  \mathcal{T} &=\textrm{add}(\begin{bmatrix}\begin{smallmatrix} P(1) \\\textrm{O}\end{smallmatrix}\end{bmatrix}_{0}\oplus\begin{bmatrix}\begin{smallmatrix} P(1) \\P(1)\end{smallmatrix}\end{bmatrix}_{1}\oplus \begin{bmatrix}\begin{smallmatrix} \textrm{O} \\P(1)\end{smallmatrix}\end{bmatrix}_{0}) \\
  \mathcal{F} &= \textrm{add}(\begin{bmatrix}\begin{smallmatrix} S(2) \\\textrm{O}\end{smallmatrix}\end{bmatrix}_{0}).
\end{align*}\par
(1) Although $i^{!}$ is exact, it is easy to see that ($\mathcal{T}$, $ \mathcal{F}$) is  a torsion pair. This illustrates that the condition of Theorem \ref{main} (1) is just a sufficient condition.\par
(2) One can compute that $i^{*}\mathcal{T}=\textrm{add}P(1)=\mathcal{T}_{1}$, $i^{!}\mathcal{F}=\textrm{add}S(2)=\mathcal{F}_{1}$, $j^{\ast}\mathcal{T}=\textrm{add}P(1)=\mathcal{T}_{1}$ and $j^{\ast}\mathcal{F}=0=\mathcal{F}_{2}$. (It illustrates Theorem
\ref{main} (2)).\par
(3) Set
\begin{align*}
  \mathcal{T} &=\textrm{add}(\begin{bmatrix}\begin{smallmatrix} P(1) \\\textrm{O}\end{smallmatrix}\end{bmatrix}_{0} ) \\
  \mathcal{F} &= \textrm{add}( \begin{bmatrix}\begin{smallmatrix} \textrm{O} \\P(1)\end{smallmatrix}\end{bmatrix}_{0}\oplus\begin{bmatrix}\begin{smallmatrix} S(2) \\\textrm{O}\end{smallmatrix}\end{bmatrix}_{0}).
\end{align*}
One can check that ($\mathcal{T}$, $ \mathcal{F}$) is  a torsion pair in $\mathcal{B}$. It is easy to check that $i_{\ast}i^{!}\mathcal{T}=i_{\ast}i^{\ast}\mathcal{T}=\textrm{add}(\begin{bmatrix}\begin{smallmatrix} P(1) \\\textrm{O}\end{smallmatrix}\end{bmatrix}_{0})=\mathcal{T}$. By Theorem \ref{main1} (1), we know that $i^{\ast}\mathcal{T}=\textrm{add}P(1)$, $i^{!}\mathcal{F}=\textrm{add}S(2)$ and ($i^{\ast}\mathcal{T}$, $i^{!}\mathcal{F}$) is a torsion pair in $\mathcal{A}$.\par
Note that $j_{!}j^{\ast}\mathcal{T}=0\subseteq \mathcal{T}$, and $j_{\ast}j^{\ast}\mathcal{F}=\textrm{add}( \begin{bmatrix}\begin{smallmatrix} \textrm{O} \\P(1)\end{smallmatrix}\end{bmatrix}_{0})\subseteq\mathcal{F}$. By Theorem \ref{main1} (2),  we know that $j^{\ast}\mathcal{T}=0$, $j^{\ast}\mathcal{F}=\textrm{add}P(1)$, and ($j^{\ast}\mathcal{T}$,$j^{\ast}\mathcal{F}$)  is a torsion pair in $\mathcal{C}$.\par
(4) Let $\mathcal{T}=\textrm{add}(\begin{bmatrix}\begin{smallmatrix} P(1) \\P(1)\end{smallmatrix}\end{bmatrix}_{1}\oplus \begin{bmatrix}\begin{smallmatrix} \textrm{O} \\P(1)\end{smallmatrix}\end{bmatrix}_{0}\oplus\begin{bmatrix}\begin{smallmatrix} S(2) \\\textrm{O}\end{smallmatrix}\end{bmatrix}_{0})$. One can check that $\mathcal{T}$ is  a rigid full subcategory of $\mathcal{B}$. For the left $\Lambda$-module $\begin{bmatrix}\begin{smallmatrix} P(1) \\\textrm{O}\end{smallmatrix}\end{bmatrix}_{0}$, there are two $\mathbb{E}$-triangles in $\mathcal{B}$
\begin{align*}
\begin{bmatrix}\begin{smallmatrix} \textrm{O} \\P(1)\end{smallmatrix}\end{bmatrix}_{0}\rightarrow\begin{bmatrix}\begin{smallmatrix} P(1) \\P(1)\end{smallmatrix}\end{bmatrix}_{1}\rightarrow \begin{bmatrix}\begin{smallmatrix} P(1) \\\textrm{O}\end{smallmatrix}\end{bmatrix}_{0}\dashrightarrow;\\
\begin{bmatrix}\begin{smallmatrix} P(1) \\\textrm{O}\end{smallmatrix}\end{bmatrix}_{0}\rightarrow\begin{bmatrix}\begin{smallmatrix} P(1) \\P(1)\end{smallmatrix}\end{bmatrix}_{1}\rightarrow \begin{bmatrix}\begin{smallmatrix} \textrm{O} \\P(1)\end{smallmatrix}\end{bmatrix}_{0}\dashrightarrow.
\end{align*}
By Remark \ref{y6}, $\mathcal{T}$ is  a cluster tilting subcategory of $\mathcal{B}$. \par
Note that $j_{\ast}j^{\ast}\mathcal{T}=\textrm{add}(\begin{bmatrix}\begin{smallmatrix} \textrm{O} \\P(1)\end{smallmatrix}\end{bmatrix}_{0})\subseteq \mathcal{T}$ and $i_{\ast}i^{\ast}\mathcal{T}=\textrm{add}(\begin{bmatrix}\begin{smallmatrix} S(2)\\\textrm{O}\end{smallmatrix}\end{bmatrix}_{0})\subseteq \mathcal{T}$. It is easy to see that $\mathcal{A}/i^{\ast}\mathcal{T}=\textrm{add}P(1)$, $\mathcal{B}/\mathcal{T}=\textrm{add}(\begin{bmatrix}\begin{smallmatrix} P(1) \\\textrm{O}\end{smallmatrix}\end{bmatrix}_{0})$, and $\mathcal{C}/j^{\ast}\mathcal{T}=0$ are abelian categories. There is a trivial recollement of abelian categories
$$
\xymatrix{\mathcal{X}\ar[rr]&&\ar@/_1pc/[ll]\ar@/^1pc/[ll]\mathcal{X}
\ar[rr]&&\ar@/_1pc/[ll]\ar@/^1pc/[ll]0}
$$
where $\mathcal{X}=\textrm{add}P(1)=\mathcal{B}/\mathcal{T}$ is an abelian category. It illustrates that the condition of Theorem
\ref{zh} is just a sufficient condition since $i^{\ast}$ is not left exact.
\end{example}

\textbf{Jian He}\\
Department of Mathematics, Nanjing University, 210093 Nanjing, Jiangsu, P. R. China\\
E-mail: \textsf{jianhe30@163.com}\\[0.2cm]
\textbf{Yonggang Hu}\\
Department of Mathematical Sciences, Tsinghua University, 100084 Beijing,  P. R. China\\
E-mail: \textsf{huyonggang@emails.bjut.edu.cn}\\[0.2cm]
\textbf{Panyue Zhou}\\
College of Mathematics, Hunan Institute of Science and Technology, 414006 Yueyang, Hunan, P. R. China.\\
E-mail: panyuezhou@163.com


\begin{thebibliography}{99}
\bibitem[BBD]{BBD} A. Beilinson, J. Bernstein, P. Deligne.
  Faisceaux pervers. (French) Analysis and topology on singular spaces, I (Luminy, 1981), 5--171, Ast\'{e}risque, 100, Soc. Math. France, Paris, 1982.

\bibitem[C]{C} J. Chen. Cotorsion pairs in a recollement of triangulated categories. Comm. Algebra 41(8): 2903--2915, 2013.

\bibitem[CZZ]{CZZ} W. Chang, P. Zhou, B. Zhu. Cluster subalgebras and cotorsion pairs in Frobenius extriangulated categories.  Algebr. Represent. Theory 22(5): 1051--1081, 2019.

\bibitem[D]{D} S. Dickson. A torsion theory for abelian categories. Trans. Amer. Math. Soc. 121: 223--235, 1966.

 \bibitem[FP]{FP} V. Franjou, T. Pirashvili. Comparison of abelian categories recollements. Doc. Math. 9: 41--56, 2004.

\bibitem[HZ]{HZ} J. He, P. Zhou. Abelian quotients of extriangulated categories. Proc. Indian Acad. Sci. 129(4): 11 pp, 2019.

\bibitem[HZZ]{HZZ} J.  Hu, D.  Zhang, P.  Zhou. Proper classes and Gorensteinness in extriangulated categories.
J. Algebra 551:  23--60, 2020.

\bibitem[IY]{IY} O. Iyama, Y. Yoshino. Mutation in triangulated categories and rigid Cohen-Macaulay modules. Inven. Math. 172: 117--168, 2008.

\bibitem[KZ]{KZ}
S. Koenig, B. Zhu.
From triangulated categories to abelian categories: cluster tilting in a general framework.
 Math. Z. 258: 143--160, 2008.


\bibitem[LW]{LW} Y. Lin, M. Wang. From recollements of triangulated categories to recollements of abelian categories. Sci. China. Math. 53(4): 1111--1116, 2010.


\bibitem[MH]{MH} X. Ma, Z. Huang. Torsion pairs in recollements of abelian categories. Front. Math. China 13(4): 875--892, 2018.

\bibitem[MV]{MV} R. MacPherson, K. Vilonen. Elementary construction of perverse sheaves. Invent. Math. 84(2):  403--435, 1986.


\bibitem[N]{N} H. Nakaoka. General heart construction on a triangulated category(1): Unifying $t$-structures and cluster tilting subcategories. Appl. Categ.Structures. 18(2): 1--21, 2009.


\bibitem[NP]{NP} H. Nakaoka, Y. Palu. Extriangulated categories, Hovey twin cotorsion pairs and model structures. Cah. Topol. G\'{e}om. Diff\'{e}r. Cat\'{e}g. 60(2): 117--193, 2019.

\bibitem[NP1]{NP1} H. Nakaoka, Y. Palu. External triangulation of the homotopy category of exact quasi-category. arXiv: 2004.02479, 2020.

\bibitem[P]{P} C. Psaroudakis. Homological theory of recollements of abelian categories. J. Algebra 398: 63--110, 2014.

\bibitem[S]{S} L. Salce. Cotorsion theories for abelian Groups. Sympos. Math. Cambridge University Press, Cambridge,  23: 11--32, 1979.

\bibitem[WWZ]{WWZ} L. Wang. J. Wei, H. Zhang. Recollements of extriangulated categories. arXiv: 2012.03258, 2020.

\bibitem[ZZ1]{ZZ1} P. Zhou, B. Zhu. Triangulated quotient categories revisited. J. Algebra 502: 196--232, 2018.

\bibitem[ZZ2]{ZZ2} P. Zhou, B. Zhu. Cluster-tilting subcategories in extriangulated categories. Theory Appl. Categ. 34: 221--242, 2019.





\end{thebibliography}
\end{document}